\documentclass[oneside,reqno]{amsart}
\usepackage{amssymb}
\usepackage{latexsym,cases}
\usepackage{paralist}

\usepackage{psfrag}
\usepackage{graphicx}
\usepackage{epsfig} 
\usepackage{float}
\usepackage{YuMacros2}
\usepackage{cite}

\usepackage{hyperref}
\hypersetup{
	colorlinks = true,
	citecolor = {red}
}

\newcommand{\PL}{\mathbb{P}_\sigma}

\numberwithin{equation}{section}

\newtheorem{thm}{Theorem}[section]

\newcommand{\lam}{\lambda_1}
\newcommand{\Om}{|\Omega|}

\newcommand{\vtil}{\tilde v}
\newcommand{\ut}{\tilde u}
\newcommand{\thetat}{\tilde\theta}

\newcommand{\ua}{u_\alpha}

\newcommand{\uL}{|u|}
\newcommand{\uuL}{|u_2|}
\newcommand{\uaL}{|u_\alpha|}
\newcommand{\ugrad}{\|u\|}
\newcommand{\uagrad}{\|u_\alpha\|}
\newcommand{\Au}{|A_0 u|}
\newcommand{\thL}{|\theta|}
\newcommand{\thgrad}{\|\theta\|}
\newcommand{\Ath}{|A_1\theta|}

\newcommand{\uLt}{|u(t)|}

\newcommand{\uaLt}{|u_\alpha(t)|}
\newcommand{\ugradt}{\|u(t)\|}

\newcommand{\thLt}{|\theta(t)|}

\newcommand{\eb}{{\bf\mathrm{e}}}

\newcommand{\ddt}[1]{\frac{\mathrm{d}{#1}}{\mathrm{d}{t}}}

\newcommand{\dd}{\;\mathrm{d}}

\newcommand{\ddz}[1]{\frac{\mathrm{d}{#1}}{\mathrm{d}{z}}}

\newcommand{\ddvp}[1]{\frac{\mathrm{d}{#1}}{\mathrm{d}{\vp}}}

\newcommand{\zmax}{z_{\max}}
\newcommand{\pmax}{q_{\max}}

\newcommand{\bbR}{\mathbb{R}}

\newcommand{\tht}{\theta}
\newcommand{\ub}{u}
\newcommand{\xb}{x}

\newcommand\vp{\vartheta}
\newcommand\vpm{\vartheta_{\max}}
\newcommand\ve{\varepsilon}
\newcommand\vt{\vartheta}
\newcommand\vtm{\vartheta_{\max}}

\newcommand\csubA{c_2}

\newcommand\cO{\mathcal{O}}

\newcommand\cA{\mathcal{A}}
\newcommand\cB{\mathcal{B}}

\newcommand\cR{\mathcal{R}}

\def\Pr{\text{Pr}}
\newcommand\Ra{\text{Ra}}
\def\Hper{H_{\text{per}}}

\begin{document}
\title[Rayleigh-B\'enard bounds]
{Algebraic Bounds on the Rayleigh-B\'enard attractor}
\author{Yu Cao$^{1}$}
\author{Michael S. Jolly$^{1,\dagger}$}
\address{$^1$Department of Mathematics\\
Indiana University\\ Bloomington, IN 47405}
\address{$\dagger$ corresponding author}
\author{Edriss S. Titi$^2$}
\address{$^2$Department of Mathematics, Texas A\&M University, 3368 TAMU,
 College Station, TX 77843-3368, USA. 
 Department of Computer Science and Applied Mathematics, Weizmann Institute
 of Science, Rehovot 76100, Israel.
 Department of Applied Mathematics and Theoretical Physics, University of Cambridge,  Cambridge CB3 0WA, UK.
} 
\author{Jared P. Whitehead$^{3}$}
\address{$^3$Department of Mathematics \\
Brigham Young University\\
Provo, UT 84602}
\email[Y. Cao]{cao20@iu.edu}
\email[M. S. Jolly]{msjolly@indiana.edu}
\email[E. S. Titi]{titi@math.tamu.edu and
 Edriss.Titi@damtp.cam.ac.uk}
 \email[J. P. Whitehead]{whitehead@mathematics.byu.edu}


\date{\today}

\begin{abstract}
The Rayleigh-B\'enard system with stress-free boundary conditions is shown to have a global attractor in each affine space where velocity has fixed spatial average.  The physical problem is
shown to be equivalent to one with periodic boundary conditions and certain symmetries. This enables a Gronwall estimate on enstrophy.  That estimate is then used to bound the $L^2$ norm of the temperature gradient on the global attractor, which, in  turn, is used to find a bounding region for the attractor in the { enstrophy--palinstrophy plane}.  All final bounds are algebraic in the viscosity and thermal diffusivity, a significant improvement over previously established estimates.  The sharpness of the bounds are tested with numerical simulations.                                                                                                                                                                                                                                                                                                                                                                                                                                                                                                                                                                                                                                                                                                                                                                                                      
\end{abstract}

\maketitle

\vspace{\baselineskip}

\normalsize

\section{Introduction}
The long-time behavior of the Rayleigh-B\'enard problem was analyzed in \cite{Foias1987attractors,temam2012infinite} for several types of boundary conditions.   In that work the authors derived explicit estimates for enstrophy and the ($L^2$-norm of the) temperature gradient on the global attractor for the case of no-slip boundary conditions in space dimension two.  They also outlined the functional setting for the case of stress-free velocity boundary conditions (see \eqref{u2zero}, \eqref{u1BC}), and mentioned that corresponding estimates can be carried out in a similar fashion.  In this paper we revisit the 2D, stress-free boundary conditions case, and as in the case of rigorous bounds on the time averaged heat transport \cite{Whitehead2011Ultimate}, we find estimates on the global attractor which are dramatically reduced from those in the no-slip boundary conditions case.  We also derive estimates for the palinstrophy and $H^2$-norm of the temperature. 

One marked difference between no-slip and stress-free boundary conditions is that in the latter case, the system is not dissipative for general initial velocity data.  This is due to the existence of steady states with arbitrarily large $L^2$-norms, namely velocity of the form $u(x,t)=(c,0)$, { (where $c$ is a constant)}
with zero temperature $\theta(x,t)=0$ such as the shear-dominated flow investigated in \cite{GoJoFlSp2014}.  Since, however, the spatial average is conserved for these flows,
the system is dissipative within each invariant affine space of fixed horizontal velocity average.  This wrinkle does not influence the estimates on the temperature or higher Sobolev norm estimates on the velocity.  

The {\it a priori} estimates are carried out in Section \ref{AprioriSect}.  The key to finding sharper bounds in the stress-free case is to extend the physical domain, as done in \cite{Farhat2017continuous}, to one that is fully periodic and twice the height of the original.  This makes the trilinear term vanish from the enstrophy balance, giving an easy bound that is $\cO(\nu^{-2})$ in terms of the kinematic viscosity $\nu$. Though the trilinear term persists when estimating the temperature gradient, we are able to avoid the exponential bound that resulted from using a uniform Gronwall lemma in \cite{Foias1987attractors}, by using the algebraic bound on the enstrophy. We find that on the global attractor
 the ($L^2$-norm of the) temperature gradient satisfies a bound that is $\cO(\Ra^2)$, for $\Pr \sim 1$, where $\Ra$ is the Rayleigh number and $\Pr$ is the Prandtl number.

We then follow the approach in \cite{Dascaliuc2010Estimates} for the Navier-Stokes equations (NSE) to obtain an estimate for the palinstrophy, with the temperature playing the role of the body force in the NSE.  This leads to curves which  bound the attractor in the {enstrophy--palinstrophy plane}, with an overall bound on palinstrophy that is $\cO(\Ra^3)$ for $\Pr \sim 1$. Using this palinstrophy bound, we then follow a similar procedure to find a bounding region for temperature $\theta$ in the { 
$\|\nabla \theta\|_{L^2}^2$--$\|\Delta \theta\|_{L^2}^2$ plane.}  

{ In Section \ref{NudgingSect}} we recall from \cite{Farhat2017continuous}  how all of these bounds impact the practicality of data assimilation by nudging with just the horizontal component of velocity of the stress-free Rayleigh-B\'enard system.  The sharpness of our rigorous bounds are tested with numerical simulations over a range of Rayleigh numbers in Section \ref{CompSect}.   Simulations are also presented there to demonstrate that the nudging algorithm works for data with much { lower} resolution than the analysis requires.  This is actually what suggested we might improve on the exponential bounds in \cite{Foias1987attractors, temam2012infinite}.  All the bounds here on the attractor are algebraic in the physical parameters.

\section{Preliminaries}\label{PrelimSect}
The Rayleigh-B\'{e}nard (RB) problem on the domain $\Omega_0=(0,L)\times(0,1)$ can be written in dimensionless form as (see, e.g., \cite{Foias1987attractors})
\begin{subequations}\label{eq-boussi0}
	\begin{gather}
	\frac{\partial u}{\partial t}
	-\nu\Delta u+(u\cdot\nabla)u+\nabla p
	=\theta \eb_2,\\
	\frac{\partial\theta}{\partial t}-\kappa\Delta\theta+(u\cdot\nabla)\theta
    ={u\cdot \eb_2},\\
	\nabla\cdot u=0,\\
	u(0;x)=u_0(x),\quad\theta(0;x)=\theta_0(x)\;,
	\end{gather}
\end{subequations}
{ where $\eb_2=(0,1)$ and $\kappa$ is the thermal diffusivity.}
In this paper, we consider the following set of boundary conditions that are stress-free on the velocity:
\begin{subequations} 
	\begin{align}
	\textrm{in the $x_2$-direction: }
	&u_2,\theta=0\ \text{at $x_2=0$ and $x_2=1$,} \label{u2zero}\\
	&\frac{\partial u_1}{\partial x_2}=0\ \text{at $x_2=0$ and $x_2=1$,} \label{u1BC}\\
	\textrm{in the $x_1$-direction: } 
	&u,\theta,p\ \text{ are of period }L\;,
	\end{align}
\end{subequations}
{ where the indices 1 and 2 refer the horizontal and vertical components, respectively.}

Following \cite{Farhat2017continuous}, in the rest of this paper, we consider the equivalent formulation of problem \eqref{eq-boussi0} subject to the fully periodic boundary conditions on the extended domain $\Omega=(0,L)\times(-1,1)$ with the following special spatial symmetries:
	\begin{align*}
	u_1(x_1,x_2)&=u_1(x_1,-x_2),\quad u_2(x_1,x_2)=-u_2(x_1,-x_2), \label{wrty} \\
	p(x_1,x_2)&=p(x_1,-x_2),\quad\quad\theta(x_1,x_2)=-\theta(x_1,-x_2)\,,
	\end{align*}
for $(x_1,x_2)\in\Omega$.
As a result of this symmetry, we observe that smooth enough functions satisfy
	\begin{align}
	u_2,\theta,\frac{\partial u_1}{\partial x_2}=0,\quad
	\textrm{for }x_2=-1,0,1\,.
	\end{align}

\subsection{Function spaces}

We will use the same notation indiscriminately for both scalar and vector Lebesgue and Sobolev spaces, which should not be a source of confusion.
We denote
\begin{align*}
(u,v):&=\int_{\Omega}u\cdot v,\quad \text{for} \ u,v\in L^2(\Omega),\\
((u,v)):&={ \sum_{i,j=1}^2\int_{\Omega} \frac{\partial u_i}{\partial x_j}\frac{\partial v_i}{\partial x_j}}, \quad \text{for} \ u,v\in H^1(\Omega),
\end{align*}
and
\[
\Norms{u}:=(u,u)^{1/2},\quad \Normd{u}:=((u,u))^{1/2}.
\]
Note that $\Normd{\cdot}$ is not a norm, { but will form part of one in \eqref{Vnorm}}.
We define function spaces corresponding to the relevant physical boundary conditions as in \cite{Farhat2017continuous}, where
\begin{quote}
	$\mathcal{F}_1$ is the set of trigonometric polynomials in $(x_1,x_2)$, with period $L$ in the $x_1$-variable, that are even, with period $2$ in the $x_2$-variable,
\end{quote}
and
\begin{quote}
	$\mathcal{F}_2$ is the set of trigonometric polynomials in $(x_1,x_2)$, with period $L$ in the $x_1$-variable, that are odd, with period $2$ in the $x_2$-variable.
\end{quote}

 The space of smooth vector-valued functions which incorporates the divergence-free condition shall be denoted by
 \begin{align*}
 	\mathcal{V}:=\{u\in\mathcal{F}_1\times\mathcal{F}_2:
 	\nabla\cdot u=0\}.
 \end{align*}
We denote the closures of $\mathcal{V}$ and $\mathcal{F}_2$ in $L^2(\Omega)$ by $H_0$ and $H_1$, respectively, which are endowed with the usual inner products
\[
(u,v)_{H_0}:=(u,v),\quad
(\psi,\phi)_{H_1}:=(\psi,\phi)
\]
and the associated norms 
\[
\Norm{u}_{H_0}:=(u,u)^{1/2},\quad 
\Norm{\psi}_{H_1}:=(\psi,\psi)^{1/2}.
\]
{ We define for $k=1,2$
\begin{align*}
\Hper^k(\Omega) =\{\phi \in H^k \ | \ \phi \ \text{has period} \ L \ \text{in}  \ x_1, \ \text{period}  \ 2  \ \text{in} \ x_2 \}\;.
\end{align*}}
Finally, we denote the closures of $\mathcal{V}$ and $\mathcal{F}_2$ in $\Hper^1(\Omega)$ by $V_0$ and $V_1$ respectively, endowed with the inner products
\[
((u,v))_{V_0}:=\frac{1}{\Om}(u,v)+((u,v)),\quad
((\psi,\phi))_{V_1}:=((\psi,\phi)),
\]
and associated norms
\begin{align}\label{Vnorm}
\Norm{u}_{V_0}:=\left(\frac{1}{\Om}\Norms{u}^2+\Normd{u}^2\right)^{1/2},\quad
\Norm{\phi}_{V_1}:=\Normd{\phi},
\end{align}
where $\Om=2L$ is the volume of $\Omega$.

\subsection{The linear operators $A_i$}
Let $D(A_0)=V_0\cap \Hper^2(\Omega)$ and $D(A_1)=V_1\cap \Hper^2(\Omega)$. Let $A_i:D(A_i)\to H_i$ ($i=0,1$) be the unbounded linear operators defined by
\begin{align*}
(A_i\phi,\psi)_{H_i}=((\phi,\psi)),\quad \phi,\psi\in D(A_i).
\end{align*}

	Due to periodic boundary conditions, we have $A_i=-\Delta$.
	The operator $A_0$ is a nonnegative operator and possesses a sequence of eigenvalues with
	\[
	0=\lambda_{0,1}\leqslant\lambda_{0,2}
	\leqslant\cdots\leqslant \lambda_{0,m}\leqslant
	\cdots,
	\]
	associated with an orthonormal basis $\{w_{0,m}\}_{m\in\N}$ of $H_0$. The operator $A_1$ is a positive self-adjoint operator and possesses a sequence of eigenvalues with 
	\[
	0<\lambda_{1,1}\leqslant\lambda_{1,2}
	\leqslant\cdots\leqslant \lambda_{1,m}\leqslant
	\cdots,
	\]
	associated with an orthonormal basis $\{w_{1,m}\}_{m\in\N}$ of $H_1$.  Observe that we have the Poincar\'{e} inequality for temperature:
	\begin{align*}
		\Norms{\theta}^2
		\leqslant \lambda_1^{-1}\Normd{\theta}^2,\quad
		&\forall \,  \theta \in V_1,\\
		\Normd{\theta}^2
		\leqslant \lambda_1^{-1}\Norms{A_1\theta}^2,\quad
		&\forall \,  \theta\in D(A_1),
	\end{align*}
where $\lambda_1=\lambda_{1,1}=\pi^2\min(1/4,L^{-2})$.

\subsection{The bilinear maps $B_i$}
Denote the dual space of $V_i$ by $V_i'$ ($i=0,1$). Define the bilinear map $B_0: V_0\times V_0\to V_0'$ (and the trilinear map $b_0: V_0\times V_0\times V_0'\to\R$) by the continuous extension of 
\[
b_0(u,v,w):=\langle B_0(u,v),w\rangle_{V_0'}:=((u\cdot\nabla)v,w), \quad
u,v,w\in\mathcal{V}.
\]
Define the scalar analogue $B_1:V_0\times V_1\to V_1'$ (and the trilinear map $b_1: V_0\times V_1\times V_1'\to\R$) by the continuous extension of 
\[
b_1(u,\theta,\phi):=\langle B_1(u,\theta),\phi\rangle_{V_1'}:=((u\cdot\nabla)\theta,\phi), \quad
u\in\mathcal{V},\ \ \theta,\phi\in\mathcal{F}_2.
\]


The bilinear maps $B_i$ (and the trilinear maps $b_i$), $i=0,1$, have the orthogonality property: 
\begin{align}\label{buvv}
	b_0(u,v,v)=0,\quad
	b_1(u,\theta,\theta)=0,\quad u,v\in V_0,\  \theta\in V_1.
\end{align}
Furthermore, due to periodicity on $\Omega$, i.e., since $A_0=-\Delta$, we have 
\begin{align}\label{buuAu}
	b_0(u,u,A_0u)=0, \quad \forall \ u\in D(A_0),
\end{align}
as well as 
\begin{align}\label{bAuuAu}
b_0(v,v,A_0w)+b_0(v,w,A_0v)+b_0(w,v,A_0v)=0\;, \quad \forall \ v,w \in D(A_0)\;,
\end{align}
(see, e.g., \cite{temam2012infinite} for \eqref{buuAu}, \cite{Foias2002statistical} for \eqref{bAuuAu}).

\subsection{Functional setting}
Following \cite{Foias1987attractors}, we have the functional form of the RB problem (\ref{eq-boussi0}):
\begin{subequations}\label{eq-benardfn0}
	\begin{gather}\label{eq-benardfn1}
	\frac{\dif u}{\dif t}+\nu A_0u+B_0(u,u)=\PL(\theta \eb_2),\\
	\label{eq-benardfn2}
	\frac{\dif \theta}{\dif t}+\kappa A_1\theta+B_1(u,\theta)
	={u\cdot \eb_2},\\
	u(0;x)=u_0(x),\quad\theta(0;x)=\theta_0(x),
	\end{gather}
\end{subequations}
where $\PL$ denotes the Leray projector.

\section{Statement of result} \label{StateSect}
\begin{thm} The Rayleigh-B\'enard problem \eqref{eq-boussi0} with stress-free boundary conditions \eqref{wrty} has a global attractor $\cA_\alpha$ within the invariant affine space 
\begin{align*}
W_\alpha=\{(u,\theta) \in V_0\times V_1: \int_\Omega u_1(x,t) \dd x =\alpha\}\;.
\end{align*}
The elements in $\cA_\alpha$ satisfy
\begin{align}\label{uLest}
\uL^2 &\le \frac{\Om}{\nu^2\lam^2} +\alpha^2\Om\;,
\end{align}
\begin{align}\label{graduest}
\ugrad^2 &\le z_{\max}:=\frac{\Om}{\nu^2\lam} \;,
\end{align}
{
\begin{align}\label{gradthetaest}
\thgrad^2 \lesssim \vartheta_{\max}:= \Om \zmax \Ra\Pr+ \left( \frac{\Om}{\lambda_1}\zmax\Ra \Pr\right)^{1/2}\;,
\end{align}
}
\begin{align}\label{p2bnd}
\Au^2 \le f(\ugrad^2)
\lesssim 
q_{\max}:=\frac{\zmax^2}{\nu^2}+\frac{\zmax^{1/2}}{\nu}\vartheta_{\max}^{1/2}\;,
\end{align}
\begin{align}\label{Athbnd}
\Ath^2 \le g(\thgrad^2) \lesssim 
\eta_{\textrm{max}}:=
\frac{z_{\max}\vartheta_{\max}}{\kappa^{2}}
+\frac{q_{\max}^{2/3}\vartheta_{\max}}{\kappa^{4/3}\lam^{1/3}}
+\frac{z_{\max}}{\kappa^{2}\lam}\;,
\end{align}
where the functions $f$, $g$ are defined below in \eqref{fdef}, \eqref{gbounds}, respectively, $\Pr$ is the Prandtl number $\nu/\kappa$, {  $\Ra=1/(\nu\kappa)$ is the Rayleigh number, and 
$Q_1\lesssim Q_2$ means $Q_1\le cQ_2$ for a nondimensional universal constant $c$ that is independent of the physical parameters.}
\end{thm}

Regions that bound the global attractor in the enstrophy--palinstrophy and $\thgrad^2$--$\Ath^2$ planes are depicted in Figures \ref{TheFig2}, \ref{etafig}, below.


\section{A priori estimates} \label{AprioriSect}
Global existence and uniqueness follows by the standard Galerkin procedure based on the trigonometric basis functions in the definitions of $\mathcal{F}_1$ and $\mathcal{F}_2$. We thus proceed with \emph{a priori} estimates.
\subsection{$L^2$ bound on temperature} 
{ We have the following maximum principle from Lemma 2.1 in \cite{Foias1987attractors}  \begin{align} \label{thetaexp}
|\theta(t)|\leq |\Omega|^{1/2}+\Theta_0e^{-\kappa t}\;,
\end{align}
where $|\Omega|$ is the volume of $\Omega$,
$$
\Theta_0 = |(\theta(0)-1)_{+}|+|(\theta(0)+1)_{-}|,
$$
and
$$
M_+=\max\{M,0\},\ M_-=\max\{-M,0\} \quad \forall \ M\in \bbR\;.
$$
While the proof in  \cite{Foias1987attractors} was done for no-slip boundary conditions,
the only place the velocity $u$ enters is the orthogonality property $b_1(u,\theta,\theta)=0$.  The proof carries over verbatim to the stress-free case by \eqref{buvv}. Consequently, we have \eqref{thetaexp} for each strong solution $(u,\theta)$ of \eqref{eq-benardfn0}.}

\subsection{$L^2$ bounds on velocity} 

We denote the space average of the horizontal velocity over the extended domain by
$$
\alpha(t)=\frac{1}{\Om}\int_\Omega \ub_1(\xb,t) \dd\xb \;.
$$
From \eqref{eq-boussi0} and the periodic boundary conditions on $\Omega$,
we find
that the spatial average of { the horizontal} velocity is conserved, i.e., $\dd \alpha/\mathrm{d} t=0$.
It follows that $\ua=u-\alpha\eb_1$ satisfies
\begin{align*}
  \ddt{\ua}+\nu A_0 \ua + B_0(\ua+\alpha\eb_1,\ua) = \theta \eb_2 \;.
\end{align*}
Since $u_\alpha$ has zero average, it satisfies the Poincar\'e inequality
\begin{align}\label{uaPoi}
	\lam |u_\alpha|^2\le
	\Normd{u_\alpha}^2\,.
\end{align}
{ Note that since $u_2$ has zero mean, it  satisfies a Poincar\'e inequality
\begin{align}\label{u2Poi}
	\lam |u_2|^2\le
	\Normd{u_2}^2\,,
\end{align}
even though $u_1$ does not.}
Taking the scalar product with $\ua$, and applying  \eqref{buvv}, the Cauchy-Schwarz and Young inequalities as well as \eqref{uaPoi}, we get  
\begin{align*}
  \frac{1}{2}\ddt{}\uaL^2+\nu\uagrad^2 
  \le \frac{1}{2\nu\lam}\thL^2 +
  \frac{\nu\lam}{2}\uaL^2 \le \frac{\thL^2}{2\nu\lam}+  \frac{\nu}{2}\uagrad^2\;.
\end{align*}
Applying \eqref{uaPoi} once again, together with \eqref{thetaexp} and Young's inequality, we have
\begin{align*}
  \ddt{}\uaL^2+\nu\lam \uaL^2
  \le \frac{1}{\nu\lam}\left(\Om +\Theta_0^2e^{-2\kappa t}\right)\;,
\end{align*}
so that 
\begin{align}\label{uaGronwall}
  \uaLt^2\le e^{-\nu\lam t}|\ua(0)|^2
  + \frac{1}{\nu\lam}\int_0^t \left(\Om +\Theta_0^2e^{-2\kappa s}\right)
   e^{\nu\lam(s-t)} \ \dd s\;,
\end{align}
and thus,
\begin{align}\label{ulim}
\limsup_{t \to \infty}\uLt^2 \le \frac{\Om}{\nu^2\lam^2} + \alpha^2\Om \;.
\end{align}


\subsection{An enstrophy bound}

We note that $\nabla u$ has zero average over $\Omega$ by the periodicity of $u$.
As a consequence, we have the Poincar\'e inequality
\begin{align}\label{poincarenabla}
  \lam \ugrad^2 \le \Au^2\;.
 \end{align}
Taking the scalar product of \eqref{eq-benardfn1} with $A_0 u$, we have by the orthogonality property \eqref{buuAu} 
\begin{align}
  \frac{1}{2} \ddt{}\ugrad^2 +\nu \Au^2 &\le |(\theta \eb_2,A_0 u)|\label{dzdt00} \\
  &\le \thL\Au
  \le \frac{1}{2\nu}\thL^2+\frac{\nu}{2}\Au^2 \;, \notag
\end{align}
hence, by \eqref{thetaexp} and \eqref{poincarenabla} we have
$$
\ddt{}\ugrad^2 +\nu\lam\ugrad^2 
\le \frac{1}{\nu}\left(\Om +\Theta_0^2e^{-2\kappa t}\right)\;,
$$
and thanks to the Gronwall inequality we obtain
\begin{align}\label{zlim}
 \limsup_{t\to \infty} \|u(t)\|^2 
 \le z_{\max}:=
  \frac{\Om}{\nu^2\lam} \;.
  \end{align}
  
Similar to the no-slip case analyzed in \cite{Foias1987attractors,temam2012infinite},  if 
$\|u(0)\| \le M_1$, $\|\theta(0)\| \le M_2$, and $\ve>0$ we have from \eqref{thetaexp}, \eqref{ulim} 
and \eqref{zlim} that there exists $t_0=t_0(M_1,M_2,\ve)$ such that
\begin{align}
\thLt^2 &\le \Om+\ve\;,  \qquad \forall \ t \ge t_0 \label{thtest}\,,\\
\uLt^2 &\le \frac{\Om}{\nu^2\lam^2} +\alpha^2\Om +\ve\;,  \qquad \forall \ t \ge t_0 \label{uest}\\
\ugradt^2 &\le \frac{\Om}{\nu^2\lam} +\ve \;, \qquad \forall \ t \ge t_0\;. \label{zest}
\end{align}
\subsection{Bound on the temperature gradient}

We start by taking the scalar product of \eqref{eq-benardfn2} with $A_1\theta=-\Delta\theta$, integrating by parts and applying the Cauchy-Schwarz and Young inequalities
\begin{align}\label{tempineq}
  \frac{1}{2} \ddt{}\thgrad^2 +\kappa\Ath^2 \le |(B_1(u,\theta),A_1 \theta)| +
  \frac{\uuL^2}{\kappa} + \frac{\kappa}{4}\Ath^2\;.
  \end{align}
We apply integration by parts to rewrite the trilinear term as
\begin{align*}
	(B_1(u,\theta),A_1 \theta)&= -\sum_{i,j=1}^{2}\int_{\Omega}u_i \partial_i\theta \partial_j^2\theta\dd x\\
&=\sum_{i,j=1}^{2}\int_{\Omega} u_i \partial_{ij}\theta \partial_j\theta  \dd x +
	\sum_{i,j=1}^{2}\int_{\Omega}\partial_j u_i \partial_i\theta \partial_j\theta\dd x\;.
	\end{align*}
We then use the chain rule to rewrite the first sum, again apply integration by parts, and then incompressibility to find
\begin{align}
\sum_{i,j=1}^{2}\int_{\Omega} u_i \partial_{ij}\theta \partial_j\theta \dd x &=
	\frac{1}{2}\sum_{i,j=1}^{2}\int_{\Omega} u_i \partial_{i}(\partial_j\theta)^2  \dd x 
	\label{incompstep} \\
&=	-\frac{1}{2}\int_{\Omega} (\partial_1u_1+\partial_2u_2)\left[(\partial_1\theta)^2+(\partial_2\theta)^2 \right] \dd x=0\;. \notag
	\end{align}
{ Applying the H\"older, Ladyzhenskaya and Young inequalities to each of the remaining four terms, we obtain}
\begin{align}
  |(B_1(u,\theta),A_1\theta)| & \le
  4\ugrad \|\nabla \theta \|_{L^4}^2  \notag \\
  & \le c_1\ugrad\thgrad\Ath \notag \\
  & \le \frac{c_1^2}{\kappa}(\ugrad\thgrad)^2 + \frac{\kappa}{4}\Ath^{2} \;.\label{B1Holder2}
\end{align}

Now combine \eqref{tempineq}, \eqref{B1Holder2} and the Poincar\'e inequality \eqref{u2Poi} so that
\begin{align*}
\ddt{}\thgrad^2 +\kappa \Ath^2 \le  \frac{2c_1^2}{\kappa}\ugrad^2\thgrad^2+
  \frac{2\uuL^2}{\kappa} \le  \frac{2c_1^2}{\kappa}\ugrad^2\thgrad^2 + \frac{2\ugrad^2}{\kappa\lam}\;.
\end{align*}
We note that by the Cauchy-Schwarz inequality and \eqref{thetaexp},
{
$$
\thgrad^2 \le \thL\Ath\le (\Om^{1/2} + \Theta_0e^{\kappa t})|A_1\theta|\;,
$$ 
}
so that 
\begin{align}\label{th4}
\ddt{}\thgrad^2 &\le  -\frac{\kappa}{(\Om^{1/2} + \ve)^2}\thgrad^4+ \frac{2c_1^2}{\kappa}\ugrad^2\thgrad^2+
  \frac{2\ugrad^2}{\kappa\lam}\;.
\end{align}
{ Let $R_1=(\Om^{1/2} + \ve)^2 $, $R_2=\zmax+\ve$, for $\ve$ as in \eqref{thtest}--\eqref{zest}.  From  \eqref{thtest}, \eqref{zest}, \eqref{th4} and Young's inequality, we have 
for all $t \ge t_0$
 \begin{align*}
\ddt{}\thgrad^2 &\le  -\frac{\kappa}{R_1}\thgrad^4+ \frac{2c_1^2}{\kappa}\thgrad^2R_2+
  \frac{2R_2}{\kappa\lam} \\
   & \le  -\frac{\kappa}{2R_1}\thgrad^4+ \frac{2c_1^4}{\kappa^3}R_1R_2^2+
  \frac{2R_2}{\kappa\lam} \\
  & \le  -\frac{\kappa}{2R_1}\left(\thgrad^4-K^4\right)\;,  
\end{align*}
where 
$$
K^4=\frac{2R_1}{\kappa}\left[\frac{2c_1^4}{\kappa^3}R_1R_2^2+
  \frac{2R_2}{\kappa\lam}\right]
$$
}
We claim that 
\begin{align}\label{claim}
\limsup_{t\to\infty} \|\theta(t)\|^2 \le 2\left[\frac{c_1^4}{\kappa^4}|\Omega|^2\zmax^2+
\frac{|\Omega|}{\kappa^2\lam} \zmax
\right]^{1/2}\;.
\end{align}
To prove this we take $\ve>0$, as above, and consider two possibilities.

Case I: \quad If $\|\theta(t)\|^2 \le (1+4\ve)^{1/2}K^2$,  for all $t \ge t_0$, then clearly
\begin{align}\label{caseclosed}
\limsup_{t\to\infty} \|\theta(t)\|^2 \le (1+4\ve)^{1/2}K^2\;, \quad \forall \ \ve > 0\;.
\end{align}

Case II: \quad Suppose there exists $t_* \ge t_0$ such that $\|\theta(t_*)\|^2 \ge (1+4\ve)^{1/2}K^2$.  We would then have that
 {
$$
\ddt{}\thgrad^2 \le  -\frac{\kappa\ve}{R_1}K^4 \;, \quad \forall \ t \ge t_*
\ \text{such that } \|\theta(t)\|^2 \ge (1+2\ve)^{1/2}K^2\;.
$$
We conclude that $\|\theta(t)\|^2$ is strictly decreasing at a rate faster than $-\kappa\ve K^4/(2R_1)$ }
for all  $t \ge t_*$
such that $\|\theta(t)\|^2 \ge (1+2\ve)^{1/2}K^2$.  In particular, there exists $t_{**}$, with 
$t_* < t_{**} < \infty$ such that $\|\theta(t_{**})\|^2 = (1+2\ve)^{1/2}K^2$.  Moreover, for all $t > t_{**}$ we have $\|\theta(t)\|^2 < (1+2\ve)^{1/2}K^2$.  As a result, we again obtain
\eqref{caseclosed}.

In either case we may now take $\ve \to 0^+$ to conclude \eqref{claim}.
Introducing the Rayleigh and Prandtl numbers in \eqref{claim} and 
{ using the concavity of the square root function,
we arrive at the bounding expression
\begin{align*}
	\limsup_{t \to \infty}\Normd{\theta(t)}^2\lesssim
	\vtm := \Om \zmax \Ra\Pr+ \left( \frac{\Om}{\lambda_1}\zmax\Ra \Pr\right)^{1/2}
\;.
\end{align*}
}
Thus, the ball $\cB_\alpha(\ve) \subset V_0\times V_1$, defined by 
\[
\cB_\alpha(\ve):=\left\{ (u,\tht) 
: \|u\|_{H^1}^2 \le \frac{1+\lam}{\nu^2\lam^2} \Om +\alpha^2\Om+ 2\ve\;,\quad
\|\theta \|^2 \le \vartheta_{\max} \right\} \;,
\] 
is absorbing.  This gives for each $\alpha$ the existence of a global attractor $\cA_\alpha$, within the invariant subspace of solutions $(u,\tht)$ where the spatial average of velocity is fixed at $\alpha$.  The global attractor is contained in $\cB_\alpha(0)$.


\subsection{Palinstrophy bound} \label{PalSect}

To estimate palinstrophy on $\cA_\alpha$ we follow  \cite{Dascaliuc2010Estimates} almost verbatim except that the effect of time independent forcing of the Navier-Stokes equations is played by the bound $\thgrad^2\le \vpm$. The other difference is that our velocity is not normalized as in \cite{Dascaliuc2010Estimates}.
For completeness, and in order to arrive at an overall bound in terms of $\nu, \kappa$, we distill the essential argument here.

Returning to \eqref{dzdt00}, we integrate by parts, and then apply the Cauchy-Schwarz inequality to get 
\begin{align*}
  -\ugrad \sqrt{\vpm} \le \frac{1}{2} \ddt{}\ugrad^2 +\nu \Au^2 \le  \ugrad\sqrt{\vpm} \;, \quad \forall\  (u,\theta) \in \cA_\alpha\;.
\end{align*}
We denote
\begin{align} \label{zpv}
  z=\|u\|^2\;, \quad 
  q=\Au^2\;, \quad 
  \zeta=|A_0^{3/2}u|^2\;,\quad 
  \vp=\|\theta\|^2\;.
\end{align}
Then whenever
\begin{align}\label{mainparabola}
  \ugrad\sqrt{\vpm} \le \frac{\nu}{2}\Au^2\;, \quad \text{equivalently} \quad
  q \ge \frac{2}{\nu}\sqrt{z\vpm}\;,
\end{align}
we have
\begin{align} \label{dzdt}
  -3\nu q \le \ddt{z} \le -\nu q\;.
\end{align}

Setting $w=A_0u$ in \eqref{bAuuAu} and applying Agmon's inequality, we have
$$
|(B_0(u,u),A_0^2u)| = |(B_0(A_0u,u),A_0u)| \le \csubA |A_0u|^2\|u\|^{1/2}
|A_0^{3/2}u|^{1/2}\;.
$$
We next take the scalar product of \eqref{eq-benardfn1} with $A_0^2u$,
and integrate by parts to obtain
\begin{align}\label{PalEq}
  \frac{1}{2} \ddt{}|A_0u|^2+\nu |A_0^{3/2}u|^2 &\le |(\theta,A_0^2u)|+
  |(B_0(u,u),A_0^2u)| \\ \notag
  & \le \|\theta\||A_0^{3/2}u| +\csubA|A_0u|^2\|u\|^{1/2}|A_0^{3/2}u|^{1/2} \;.
  \end{align}
Note that by the Cauchy-Schwarz inequality
\begin{align}\label{zetavp0}
  \zeta:=|A_0^{3/2}u|^2 \ge \frac{\Au^4}{\|u\|^2} = \frac{q^2}{z}
  \ge \frac{4}{\nu^2}\vpm
\end{align}
in the region
\begin{align}\label{Rregion}
  \cR:=\left\{ (z,q) \ \big| \ q \ge \frac{2}{\nu} \sqrt{z\vpm}\ \right\}\;.
\end{align}
It follows that
\begin{align*}
  \|\theta\||A_0^{3/2}u| = \vp^{1/2}\zeta^{1/2}\le \frac{\nu}{2} \zeta \qquad
  \forall \ (z,q) \in \cR \;,
\end{align*}
and hence, as in \cite{Dascaliuc2010Estimates},
\begin{align}\label{dpdt}
  \ddt{q} \le \psi(\zeta):=-\nu \zeta +2\csubA qz^{1/4}\zeta^{1/4} \;.
  \end{align}
To close the system (eliminate $\zeta$) we find that the maximum of $\psi$
is achieved at
\begin{align*}
  \zeta_{\text{max}}:=\left(\frac{\csubA}{2\nu}qz^{1/4}\right)^{4/3} \quad \text{with a value} \quad
  \psi_{\text{max}}=3\nu\zeta_{\text{max}}\;.
\end{align*}
We note that
\begin{align*}
  \frac{q^2}{z} \ge \zeta_{\text{max}} \quad \text{if and only if} \quad
  q \ge \left(\frac{\csubA}{2\nu}z\right)^2
\end{align*}
so that by \eqref{zetavp0}
 \begin{align}\label{dpdt1}
   \ddt{q} \le
     \psi_{\text{max}}= \frac{3}{\nu^{1/3}}\left(\frac{\csubA}{2}qz^{1/4}\right)^{4/3} 
    \quad \text{if} \quad q \le \left(\frac{\csubA}{2\nu} z \right)^2
 \end{align}
 and
 \begin{align}\label{dpdt2}
   \ddt{q} \le \psi(q^2/z)=-\nu\frac{q^2}{z}+2\csubA q^{3/2} &
    \quad \text{if} \quad q \ge \left(\frac{\csubA}{2\nu} z \right)^2 \;.
   \end{align}
 We see that
 \begin{align}\label{pdecr}
   \ddt{q} \le 0 \quad \text{if}  \quad q \ge \left(\frac{2\csubA}{\nu} z\right)^2
   \quad \text{and} \quad q \ge \frac{2}{\nu}\sqrt{z\vpm}\;.
   \end{align}

\begin{figure}[h]
\psfrag{z0}{$\frac{\Om}{\nu^2\lam}$}
\psfrag{z}{$z=\ugrad^2$}
\psfrag{z1}{$z_1$}
\psfrag{z2}{$z_2$}
\psfrag{root}{$\frac{2}{\nu}\sqrt{z\vpm}$}
\psfrag{th}{$\thgrad^2$}
  \psfrag{sqrt}{$\frac{2(\Om)^{1/2} }{\kappa\lam^{1/2}}\ugrad$}
  \psfrag{line}{$\frac{8c L}{\kappa^2}\ugrad^2$}
  \psfrag{b2}{$\frac{32cL^2}{\kappa^2\nu^2\lam}$}
  \psfrag{b1}{$ $}
  \psfrag{k}{$\frac{\kappa^{2}}{\Om}$}
  \psfrag{par1}{$\left(\frac{c_2z}{2\nu}\right)^2$}
  \psfrag{par2}{$\left(\frac{2c_2z}{\nu}\right)^2$}
  \psfrag{Poincare}{$\lam z$}
\psfrag{z1}{$z_1$}
\psfrag{z2}{$z_2$}  
\psfrag{p0}{$q_0$}
  \psfrag{p}{$q=\Au^2$}
  \psfrag{pbnd}{$q_2$}
  \psfrag{f1}{$f_1$}
  \psfrag{f2}{$f_2$}
  \psfrag{f3}{$f_3$}
  \centerline{  \includegraphics[scale=.65]{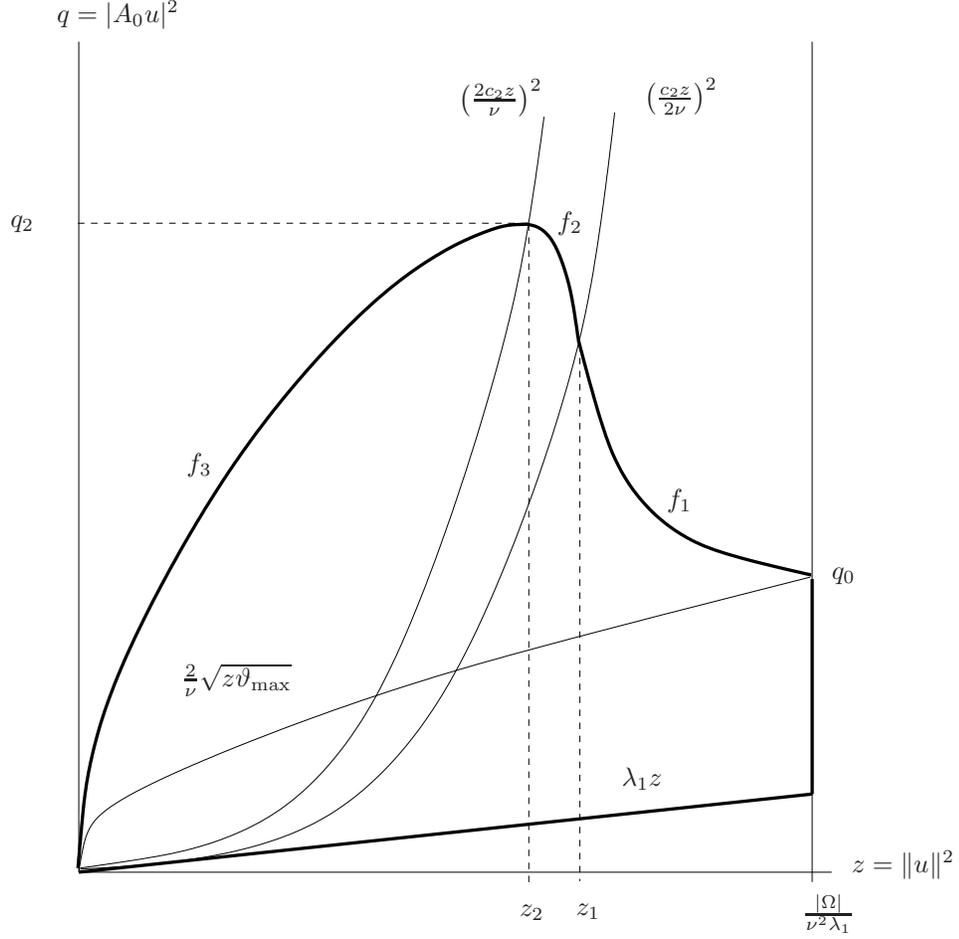}}
\caption{Qualitative sketches of the curves bounding $\cA_\alpha$.}
\label{TheFig2}
\end{figure}

 By considering the steepest descent possible below 
 $$
 q =\left(\frac{2\csubA}{\nu} z\right)^2
 $$
 and the most shallow ascent possible above this parabola, we find three bounding curves
 $q=f_j(z)$, $j=1,2,3$, after solving, in order, three final value problems.
 The first combines the (positive) bound in \eqref{dpdt1} with the upper bound in
 \eqref{dzdt}
 \begin{align*}
&\ddz{q} =-3\left(\frac{\csubA}{2\nu}\right)^{4/3}(qz)^{1/3}\;,  \qquad \text{for} \quad
     z_1\le z \le z_0=\zmax =\frac{|\Omega|}{\nu^2\lam}\\
     &q(z_0)=q_0:=\frac{2}{\nu^2}\zmax^{1/2}\vpm^{1/2}\;.
 \end{align*}
  The second picks up where the first leaves off and combines the (positive) bound in \eqref{dpdt2} with the upper bound in \eqref{dzdt}
  \begin{align*}
  & \ddz{q} =\frac{q}{z} -\frac{2\csubA}{\nu}q^{1/2}\;,  \qquad \text{for} \quad
     z_2\le z \le z_1 \\
     & q(z_1)=q_1\;, 
  \end{align*}
  while the third combines the (negative) bound in \eqref{dpdt2} with the {\it lower}
  bound in \eqref{dzdt}
  \begin{align*}
 & \ddz{q} =\frac{q}{3z} -\frac{2\csubA}{3\nu}q^{1/2}\;,  \qquad \text{for} \quad
     0\le z \le z_2 \\
      &q(z_2)=q_2\;, 
  \end{align*}
  where $q_1$, $q_2$ are determined by the intersections of $f_1$ and $f_2$ (defined below) with the parabolas
  $$
  q=\left(\frac{\csubA}{2\nu} z\right)^2 \quad \text{and} \quad  q=\left(\frac{2\csubA}{\nu} z\right)^2
  $$
  respectively.
This results in a convex function in $z$ 
\begin{align*} 
f_1(z):=\left[\frac{3}{2}\left(\frac{\csubA}{2\nu}\right)^{4/3}\left(z_0^{4/3}-z^{4/3}\right) + q_0^{2/3} \right]^{3/2}
\end{align*}
and concave functions in $z$
\begin{align*} 
f_2(z)&:=\frac{1}{\nu^2}\left[-2\csubA z + \left(\nu q_1^{1/2}+2\csubA z_1\right)
\left(\frac{z}{z_1}\right)^{1/2}\right]^{2}\;, \\
f_3(z)&:=\frac{1}{25\nu^2}\left[-6\csubA z + \left(5\nu q_2^{1/2}+6\csubA z_2\right)
\left(\frac{z}{z_2}\right)^{1/6}\right]^{2} \;.
\end{align*}
A qualitative sketch of these three curves is shown in Figure \ref{TheFig2}.  It is shown in \cite{Dascaliuc2010Estimates} that the curve $q=f_3(z)$ does not intersect the curve $q=2\sqrt{z\vpm}/\nu$.
Let 
\begin{align} \label{fdef}
f(z):=\begin{cases} f_1(z) & \text{if} \quad z_1 \le z \le \zmax \\
   f_2(z) & \text{if} \quad z_2 \le z < z_1 \\
   f_3(z) & \text{if} \quad 0 \le z \le z_2 \;.
   \end{cases}
   \end{align}

To prove \eqref{p2bnd}, suppose there is an element in $\cA_\alpha$ such that $q(0) > f(z(0))$. The solution through any element in $\cA_\alpha$ exists for all negative time.
If $q(t)>f(z(t))$ for all $t<0$, since $q(t)$ increases with negative time, as long as $z(t)<z_2$, we have $q(t)>\min\{q(0),q_0\}$.  By the upper bound in \eqref{dzdt},  $z(t)$ would then exceed $z_2$ in finite negative time.  Thus, we must have $q(t) \le f(z(t)) $ at some $t<0$. But forward in time, the region $q \le f(z)$ is invariant, contradicting the assumption that the initial condition satisfied  $q(0) > f(z(0))$.
  
We now find an overall bound on palinstrophy in $\cA_\alpha$. 
A straightforward calculation shows that substituting 
$$
q_1=\left(\frac{c_2z_1}{2\nu}\right)^2 \quad \text{into} \quad q_2= f_2(z_2)=\left(\frac{2c_2z_2}{\nu}\right)^2
$$
reduces to
  \begin{align*}
    z_2=\frac{25}{64}z_1\;.
\end{align*}
Similarly, using
$$
q_0=\frac{2}{\nu}\sqrt{z_0\vpm} \quad \text{in} \quad q_1=f_1(z_1)=\left(\frac{c_2z_1}{2\nu}\right)^2 \;,
$$
leads to 
$$
z_1\le\left[\frac{3}{2}z_0^{4/3}+ \frac{4\nu^{2/3}}{c_2^{4/3}}z_0^{1/3}\vpm^{1/3}\right]^{3/4} 
\lesssim 
\zmax+\nu^{1/2}\zmax^{1/4}\vartheta_{\max}^{1/4}$$
so that
\begin{align}\label{p22bnd}
\Au^2 \le q_2 \lesssim 
q_{\max}:=\frac{\zmax^2}{\nu^2}+\frac{\zmax^{1/2}}{\nu}\vartheta_{\max}^{1/2}
\;.
\end{align}

\subsection{A bound on $\Ath$}
From \eqref{tempineq} and \eqref{B1Holder2} we have 
\begin{align*}
-\frac{c_3}{\kappa}\ugrad^2\thgrad^2
  -\frac{2\ugrad^2}{\kappa\lam} \le   \ddt{}\thgrad^2 +\kappa \Ath^2 \le  \frac{c_3}{\kappa}\ugrad^2\thgrad^2+
  \frac{2\ugrad^2}{\kappa\lam}\;.
  \end{align*}
Thus, if 
\begin{align*}
    \Ath^2 \ge \frac{2}{\kappa^2}\left( c_3\ugrad^2\thgrad^2+
  \frac{2\ugrad^2}{\lam}\right)\;, 
  \end{align*}
  it follows that
  \begin{align}\label{dAth}
-\frac{3}{2}\kappa \Ath^2 \le    \ddt{}\thgrad^2 \le -\frac{1}{2}\kappa \Ath^2  \;.
\end{align}
  
We next take the scalar product of the temperature equation with $A_1^2\theta=\Delta^2\theta$ and using the fact that $A_0u=-\Delta u$, write
  \begin{align}\label{dAth1}
  \frac{1}{2}\ddt{}\Ath^2 + \kappa |A_1^{3/2} \theta|^2 &\le \Au\Ath  
  + |(B_1(u,\theta),A_1^2\theta)| \;.
  \end{align}
We need to move two derivatives in the trilinear term in order to ultimately obtain a bound for it in which the highest order norm is $|A_1^{3/2}\theta|$.  We integrate by parts to write
\begin{align*}
	(B_1(u,\theta),A_1^2 \theta)&= \sum_{i,j,k=1}^{2}\int_{\Omega}u_i \partial_i\theta \partial_j^2\partial_k^2\theta\dd x\\
&=-\sum_{i,j,k=1}^{2}\int_{\Omega} u_i \partial_{ij}\theta \partial_j\partial_k^2\theta  \dd x -
	\sum_{i,j,k=1}^{2}\int_{\Omega}\partial_j u_i \partial_i\theta \partial_j\partial_k^2\theta\dd x =I+II \;.\\
	\end{align*}
We then integrate the first summation by parts	
	\begin{align*}
I=\sum_{i,j,k=1}^{2}\int_{\Omega} u_i \partial_{i}\partial_j^2\theta \partial_k^2\theta \dd x 
+ \sum_{i,j,k=1}^{2}\int_{\Omega} \partial_j u_i \partial_{i,j}\theta \partial_k^2\theta \dd x =I_a+I_b
	\end{align*}
and split the resulting first summation as	
	\begin{align*}
I_a=\sum_{i,j=1}^{2}\int_{\Omega} u_i \partial_{i}\partial_j^2\theta \partial_j^2\theta \dd x 
+ \sum_{i,j\neq k=1}^{2}\int_{\Omega} u_i \partial_{i}\partial_j^2\theta\partial_k^2\theta \dd x  =I_{a_1}+I_{a_2}
	\end{align*}
Proceeding as in \eqref{incompstep}, we find that $I_{a_1}=0$.  Integrating by parts again, we have	
\begin{align*}
I_{a_2}=-\sum_{i,j \neq k=1}^{2}\int_{\Omega} \partial_{i}u_i \partial_j^2\theta \partial_k^2\theta \dd x -
	\sum_{i,j\neq k=1}^{2}\int_{\Omega} u_i \partial_j^2\theta  \partial_{i}\partial_k^2\theta  \dd x\;.
	\end{align*}
	Since the first sum is zero by incompressibility,  we have by symmetry that $I_{a_2}=-I_{a_2}$, and thus $I_{a_2}=0$.   Integrating by parts one more time, 
	we have	
\begin{align*}
II=\sum_{i,j,k=1}^2\int_{\Omega}  \partial_j^2u_i\partial_i\theta \partial_k^2\theta \dd x 
+ I_b	\;.
	\end{align*}
	After gathering what remains, we use Agmon's and Ladyzhenskaya's inequalities to estimate the trilinear term as
\begin{align*}
|(B_1(u,\theta),A_1^2\theta)|&=	|\sum_{i,j,k=1}^2\int_{\Omega}  \partial_j^2u_i\partial_i\theta \partial_k^2\theta \dd x  + 
	2\sum_{i,j,k=1}^{2}\int_{\Omega} \partial_j u_i \partial_{i,j}\theta \partial_k^2\theta \dd x | \\
	&\le c\Au \thgrad^{1/2} |A_1^{3/2}\theta|^{1/2}\Ath + 
	   c\ugrad^{1/2}\Au^{1/2}\|\theta\|_{H^2}\Ath^{1/2}|A_1^{3/2}\theta|^{1/2} \\
	&\le c_4{\Au}\frac{\Ath^{3/2}}{\lam^{1/4}}|A_1^{3/2}\theta|^{1/2}
	   = \frac{c_4}{\lam^{1/4}}q^{1/2}\eta^{3/4}\xi^{1/4}\;,
	\end{align*}
	where 
	$\eta=\Ath^2$, $\xi= |A_1^{3/2}\theta|^2$ and for convenience in what follows, we take  $c_4=2\max(c,c_3)$.

  Using this in \eqref{dAth1}, we find
  \begin{align*}
  \frac{1}{2}\ddt{}\Ath^2 + \kappa |A_1^{3/2} \theta|^2 &\le \Au\Ath + 
  \frac{c_4}{\lam^{1/4}} \Au\Ath^{3/2}|A_1^{3/2}\theta|^{1/2} \\
  &\le 
  \frac{2c_4}{\lam^{1/4}} \Au\Ath^{3/2}|A_1^{3/2}\theta|^{1/2} \;.
  \end{align*}
 Thus, invoking our palinstrophy bound $\pmax$, we have
  \begin{align*}
  \ddt{}\eta \le \Phi(\xi) :=-2\kappa \xi +\frac{4c_4}{\lam^{1/4}}\pmax^{1/2}\eta^{3/4}\xi^{1/4}\;.
  \end{align*}
  We find that 
  $$
  \Phi(\xi) \le \Phi_{\text{max}}=\frac{2}{\kappa^{1/3}}\left(\frac{c_4}{2\lam^{1/4}}\right)^{4/3} \pmax^{2/3} \eta
  $$
  and that 
 $$
  \Phi(\xi) \le 0 \quad \forall \ \xi \ge \xi^*:=\gamma\eta \;,
\quad
 \text{where}\quad \gamma:=\left(\frac{2c_4}{\kappa\lam^{1/4}}\right)^{4/3} \pmax^{2/3}\;.
   $$ 
In terms of $z_0$, our enstrophy bound on the attractor, \eqref{dAth} holds for
\begin{align}\label{A1thcond}
\eta \ge g_3(\vp):=\frac{z_{\max}}{\kappa^2}\left( c_4\vt+
  \frac{4}{\lam}\right) \;.
  \end{align}

 Once again, by the Cauchy-Schwarz inequality, we have 
 $$
 |A_1^{3/2}\theta| \ge \frac{\Ath^2}{\thgrad}\;, \quad \text{i.e.,} \quad
\xi \ge  \frac{\eta^2}{\vt}\;.
$$
Thus for 
$$
 \frac{\eta^2}{\vt} \le \xi^*\;, \quad \text{equivalently} \quad \eta \le \gamma \vt\;, 
 $$
 we combine
 \begin{align} \label{detadvp}
 \ddt{}\eta \le \Phi_{\text{max}}\quad 
\text{with}\quad 
 \ddt{}\vt\le -\frac{\kappa}{2} \eta
\end{align} 
and solve 
$$
\ddvp{\eta} =-\gamma_0\;, \quad \eta(\vtm)=\eta_0:=\frac{\zmax}{\kappa^2}\left( c_4\vtm+
  \frac{4}{\lam}\right)\;, \quad \text{where}\quad \gamma_0=4^{-1/3}\gamma 
$$
to find a straight-line solution
$$
\eta=g_1(\vt):=\eta_0-\gamma_0(\vt-\vtm) \;.
$$
We then find the intersection of this line with $\eta=\gamma\vt$ to be at $(\vt_1,\eta_1)$,
where
\begin{align}\label{eta1}
\vt_1=\frac{c_4 \zmax/\kappa^2 + \gamma_0}{\gamma+\gamma_0} \vtm + 
\frac{4 \zmax}{\kappa^2\lam(\gamma+\gamma_0)}\;, \quad \text{and} \quad
\eta_1=\gamma \vt_1\;.
\end{align}

For $\eta \ge \gamma \vt$ we combine
$$
\ddt{}\eta \le \Phi(\eta^2/\vp)=-2\kappa \frac{\eta^2}{\vp} 
     + \frac{4c_4}{\lam^{1/4}}\pmax^{1/2} \frac{\eta^{5/4}}{\vp^{1/4}} \quad \text{with}
     \quad 
     \ddt{}\vp \ge -\frac{3}{2} \kappa \eta
     $$
     and solve 
     \begin{align*}
     \ddvp{\eta}= \frac{4}{3\vp} \eta - \frac{8c_4}{3\lam^{1/4}\kappa} \pmax^{1/2}\frac{\eta^{1/4}}{\vp^{1/4}}
     \end{align*}
     to find
     \begin{align*}
     \eta =g_2(\vp):=\left[\left(\frac{\vp}{\vp_1}\right)\eta_1^{1/4} +\tilde\gamma
     \left(\vp^{3/4}- \vp\vp_1^{-1/4}\right)\right]^{4/3}\;,
     \end{align*}
     where 
     $$
     \tilde\gamma=\frac{8c_4}{\lam^{1/4}\kappa} \pmax^{1/2}\;.
     $$
     
\begin{figure}[h]
\psfrag{vpm}{$\vpm$}
\psfrag{vp}{$\vp=\thgrad^2$}
\psfrag{vp1}{$\vp_1$}
\psfrag{Poincare}{$\lam\vp$}
  \psfrag{eta}{$\eta=\Ath^2$}
  \psfrag{eta1}{$\eta_{\text{max}}$}
   \psfrag{gamma}{$\gamma\vp$}
  \psfrag{g1}{$g_1$}
  \psfrag{g2}{$g_2$}
  \psfrag{g3}{$g_3$}
  \centerline{\includegraphics[scale=.5]{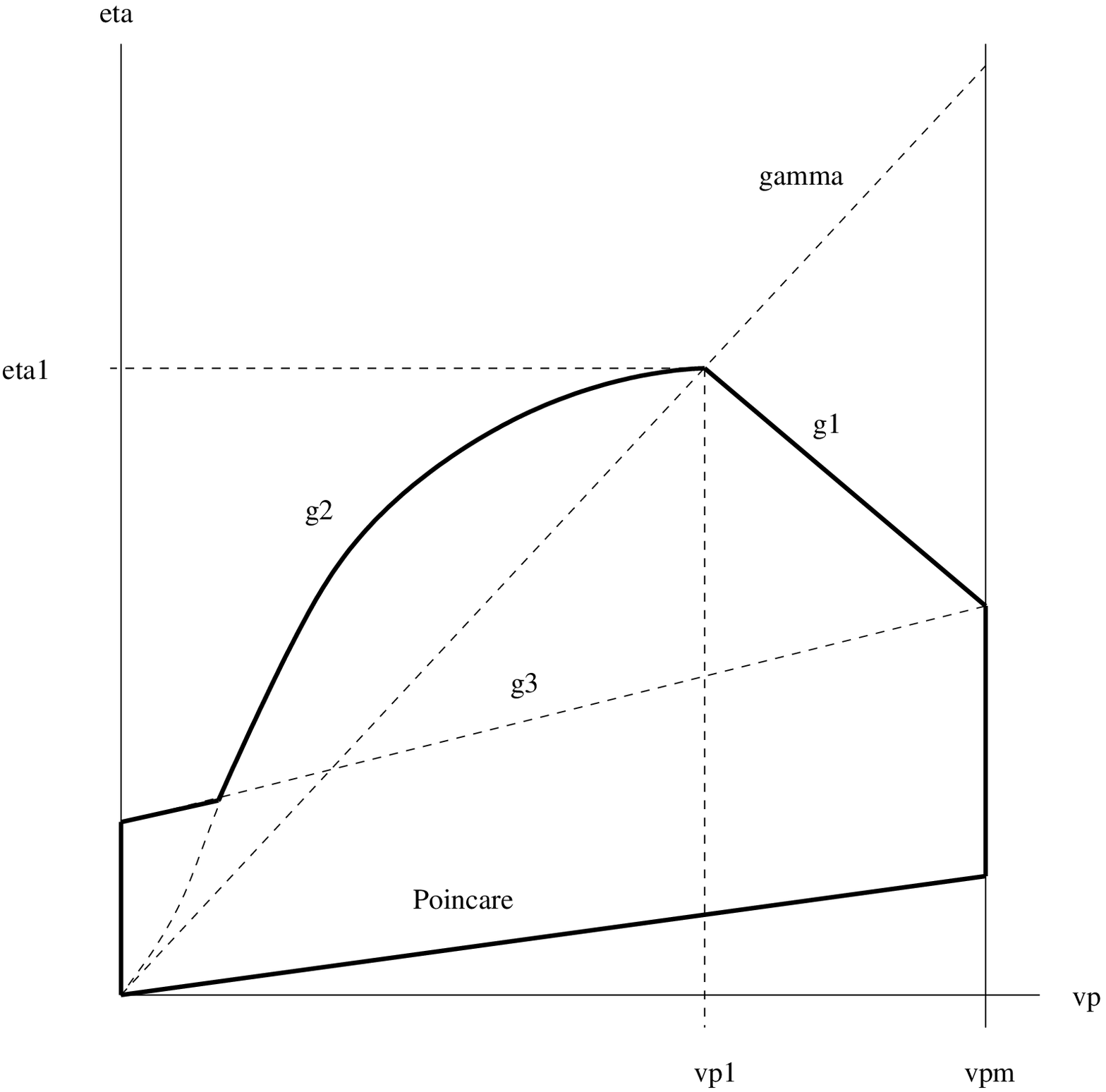}}
\caption{Bounding region in the { $\thgrad^2$--$\Ath^2$ plane}.}
\label{etafig}
\end{figure}

As we argued in Section \ref{PalSect}, if an element in the global attractor were to project in the { $\vp$--$\eta$ plane} above 
\begin{align}\label{gbounds}
\eta=\max\left\{g_1(\vp),g_2(\vp),g_3(\vp)\right\}\;,
\end{align}
then by \eqref{detadvp} the solution through it would, in finite negative time, have to enter the region below the curves in \eqref{gbounds}.  Yet, this region is invariant.  We conclude
from \eqref{eta1} and \eqref{p22bnd} that we have an overall bound on the global attractor of 
\begin{align*}
\Ath^2 \le \eta_1 \lesssim \eta_{\max} :=
\frac{z_{\max}}{\kappa^2}\vartheta_{\max}
+\gamma\vartheta_{\max}
+\frac{\zmax}{\kappa^2\lam}\;.
\end{align*}
A qualitative sketch of the region bounding the global attractor in this plane is shown in Figure \ref{etafig}.

\section{Implications for data assimilation}  \label{NudgingSect}
Suppose reality is represented by a particular solution to an evolution equation 
 \begin{equation*}\label{evol}
\frac{dv}{dt}  
    = F(v)  \;,
\end{equation*}
where the initial data $v(0)$ is \textit{not} known.  Instead continuous data of the form $I_hv(t)$ is known over an interval, $t\in [t_1,t_2]$, for a certain type of interpolating operator $I_h$ with spatial resolution $h$.
The nudging approach to data assimilation  amounts to solving the auxiliary  system
 \begin{align}\label{feedback:dissip}
        \frac{d\vtil}{dt}=F(\vtil)-\mu I_h(\vtil-v)\;,
 \end{align}
 using \textit{any} initial condition, e.g., $\vtil_0=0$.   It was shown in \cite{Azouani2013feedback, Azouani2014continuous} that if $\mu>0$ is sufficiently large, and correspondingly, $h$ sufficiently small, then  $v(t)-\vtil(t) \to 0$, in some norm, at an exponential rate, as $t \to \infty$.   
 In fact, computations indicate that this approach works with data that is much more coarse than suggested by rigorous estimates (see \cite{Altaf2017Down,Farhat2018Assim,Gesho2016,diLeoni2018Inferring}).
 Flexibility in the choice of interpolant is one of the main advantages of injecting the observed data through a feedback nudging term, rather than into terms involving spatial derivatives \cite{Azouani2014continuous, Jones1992Determining}. 
 Numerical errors are shown to be bounded uniformly in time for semi-discrete \cite{Mondaini2018Uniform} and fully discrete schemes \cite{Ibdah2019Uniform} for \eqref{feedback:dissip}.

Now consider this approach for the stress-free Rayleigh-B\'enard system
\eqref{eq-benardfn0} using data from only the horizontal component of velocity.  This means solving the auxiliary system
	\begin{gather*}
	\frac{\dif \ut}{\dif t}+\nu A_0\ut+B_0(\ut,\ut)=\PL(\thetat \eb_2)-\mu \PL I_h(\ut_1-u_1)\eb_1,\\
	\frac{\dif \thetat}{\dif t}+\kappa A_1\thetat+B_1(\ut,\thetat)
	={\ut\cdot \eb_2},\\
	\ut(0;x)=0,\quad\thetat(0;x)=0.
	\end{gather*}
It was proved in \cite{Farhat2017continuous} that if $\mu h^2 \lesssim \nu$ and 
\begin{align}\label{condAu}
\mu \ge K_1 \sim \frac{1}{\kappa\lam} + \frac{1}{\nu\kappa^2} + \frac{1}{\kappa} +
\frac{\Au^2}{\nu}\;,
\end{align} 
then
$$
\|u(t) - \tilde u(t)\| + |\theta(t)-\tilde\theta(t)| \to 0 \quad \text{as} \quad t \to \infty
$$
at an exponential rate.
Also shown there was that if 
\begin{align}\label{condAth}
\mu \ge K_2 \sim K_1 +\frac{1}{\kappa}\thgrad^2\Ath^2\;,
\end{align}
then the stronger convergence
$$
\|u(t) - \tilde u(t)\| + \|\theta(t)-\tilde\theta(t)\| \to 0, \quad \text{as} \quad t \to \infty,
$$
holds at an exponential rate.  The bounds in this paper on $\|\theta\|$, $|A_0u|$ and $|A_1\theta|$ are all algebraic, suggesting that data assimilation by nudging with just the horizontal velocity could be effective for the stress-free Rayleigh-B\'enard system.  We present computational evidence to this effect in the next section.

\section{Computational Results}\label{CompSect}

The computations presented below were done using Dedalus,  an open-source package for solving partial differential equations using pseudo-spectral methods (see \cite{dedalus}). The time stepping is done by a four-stage third order Runge-Kutta method.

We solve  \eqref{eq-boussi0} with $L=2$ in the physical domain $\Omega_0=(0,L)\times (0,1)$.  The physical parameters of viscosity and thermal diffusivity are related to the Rayleigh and Prandtl numbers through
\[
\nu = \sqrt{\frac{\Pr}{\Ra}},\quad \kappa = \frac{1}{\sqrt{\Ra\cdot \Pr}}\;.
\]
We take $\Pr=1$ so that in our dimensionless variables $\Ra:=(\nu\kappa)^{-1}=\nu^{-2}$ and use $n_F$ Fourier modes in the $x_1$-direction and $n_C$ Chebyshev modes in the $x_2$-direction. The numbers of modes used are $n_F\times n_C=256\times128$, $1024\times512$, and $2048\times1024$ for runs at $\Ra=10^6, 10^7, 10^8$ respectively. 

\subsection{Sharpness} \label{sharpness}

\begin{figure}[!htbp] 
	\includegraphics[width=\textwidth]{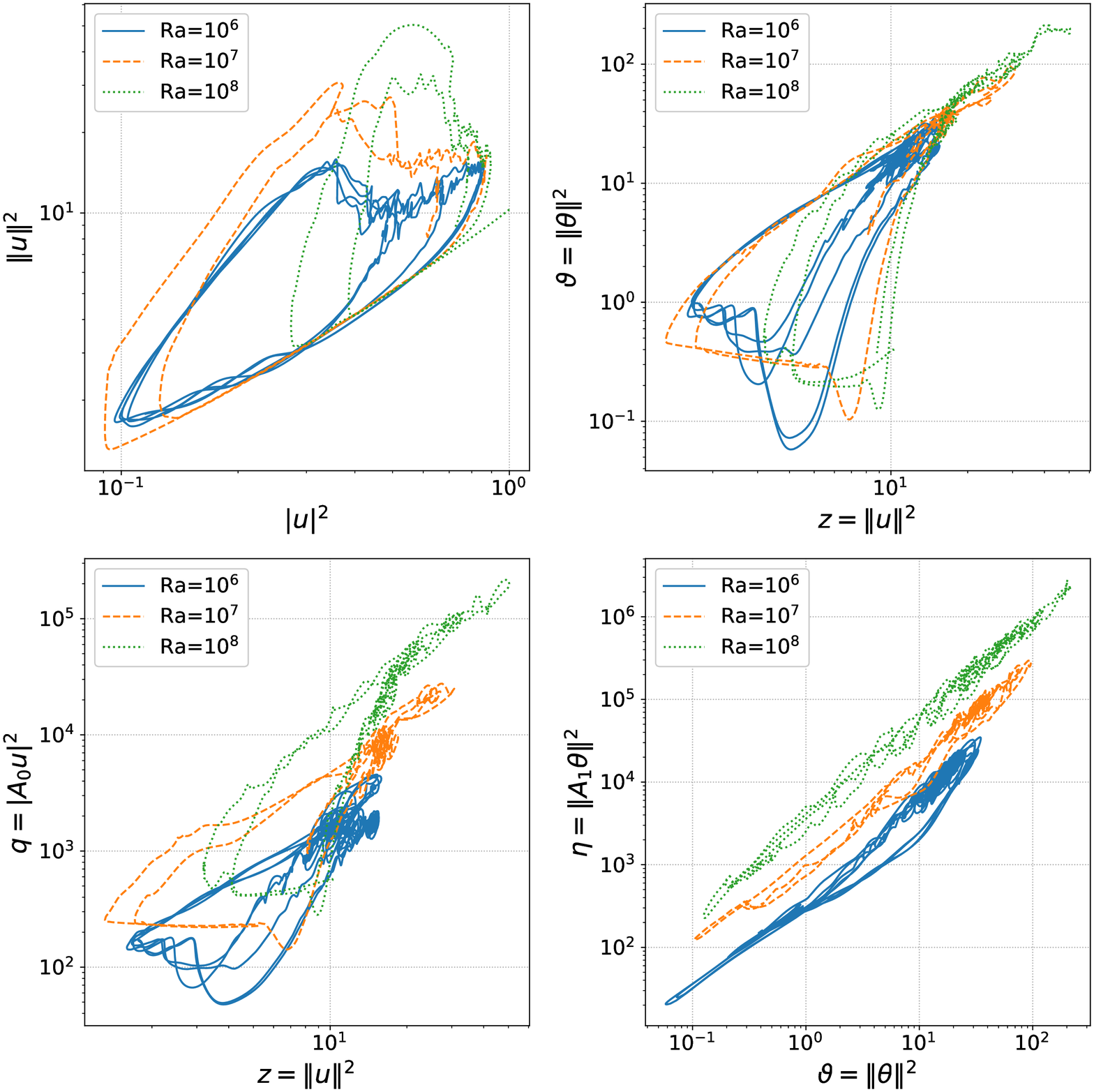}
	\caption{Projections after a transient period}
	\label{attr}
\end{figure}

\begin{figure}[!htbp] 
	\includegraphics[width=\textwidth]
	{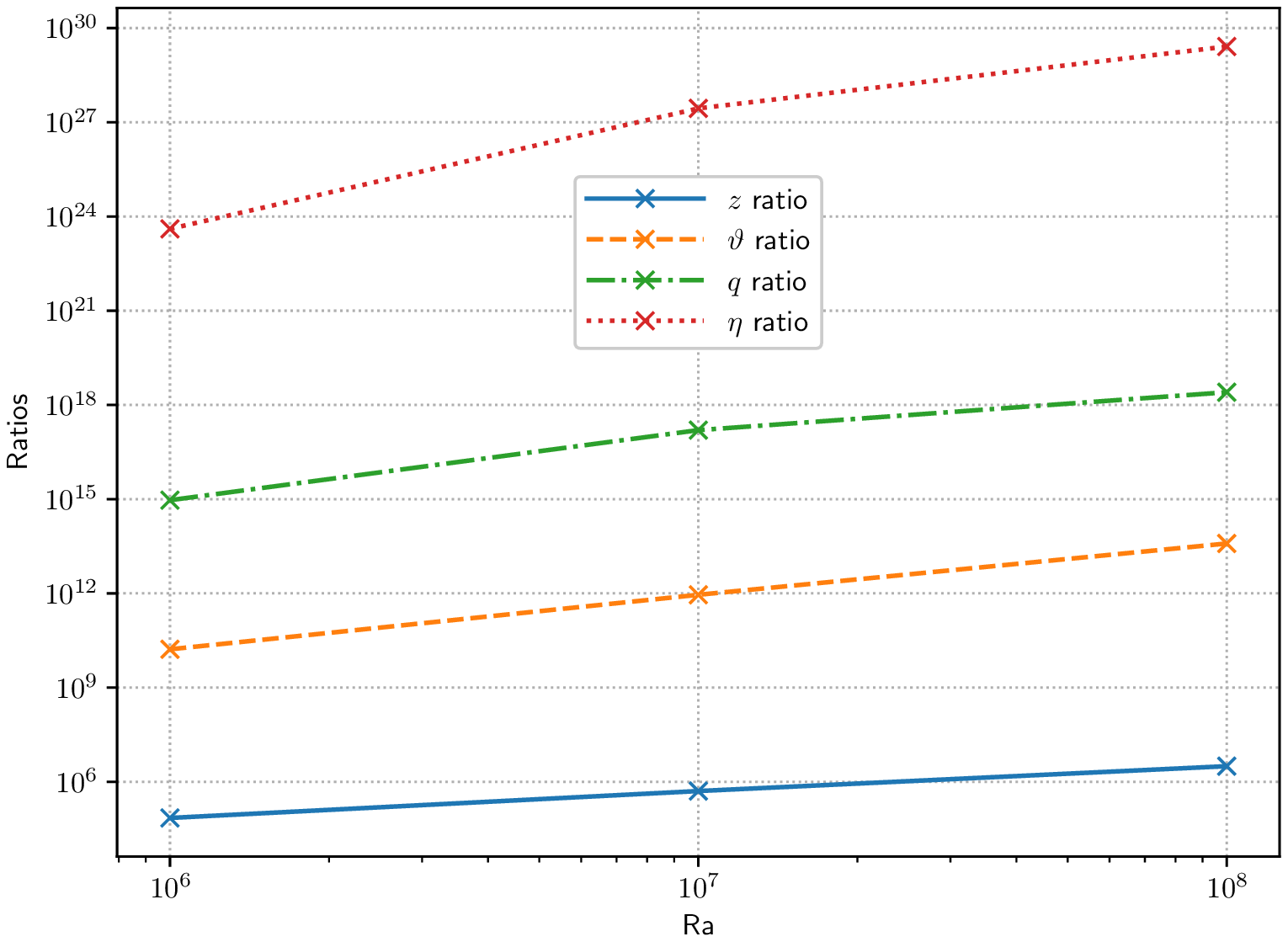}
	\caption{Ratios in \eqref{rats}}
	\label{ratios}
\end{figure}


Each plot in Figure \ref{attr} shows the projection of a solution after a transient phase in a plane spanned by the norms bounded by our analysis.  The solutions are plotted over the time period $200\le t\le1000$ for $\Ra=10^6, 10^7$ and over $200\le t\le1485$ for $\Ra=10^8$ (time units in the RB system \eqref{eq-boussi0}).  The initial condition in each case is $(u_0,\theta_0)=(0,0)$ so the average $\alpha$ of the horizontal velocity is zero.

It is not surprising that our rigorous overall bounds as well as the curves in Figures \ref{TheFig2}, \ref{etafig} are orders of magnitude greater than the norms of these solutions.  Plotting the bounds and curves together with the solutions is not revealing.  Instead, to see a trend in sharpness, we plot in Figure \ref{ratios} the ratios
\begin{align}\label{rats}
\frac{\zmax}{\max_{\cA} z}\;, \qquad \frac{\vt_{\max}}{\max_{\cA} \vt}\;, \qquad \frac{\pmax}{\max_{\cA} q}\;, \qquad \frac{\eta_{\max}}{\max_{\cA} \eta}\;.
\end{align}
{ Using the numerical values for the $z$ ratio, we gauge the (highest) power in  \eqref{graduest} to be inflated (at least over this range of the Rayleigh number),  by an addition of $\beta$, where 
$$
100^\beta = \frac{3.14\times 10^6}{7.05\times 10^4}\;, \quad \text{i.e.,} \ \beta=.824\;.
$$ 
A similar calculation for the $\vt$ ratio gives $\beta=1.68$.  We note that the curves are bending favorably for the ratios for $q=\text{palinstrophy}$ and $\eta=\Ath$.}
 
\subsection{Data assimilation} 
Nudging is carried out at $\Ra=10^6$ using the interpolant operator $I_h$ at every $m^{\rm{th}}$ nodal value in each direction, i.e., 
\[
h(m)= \max\{h_F(m),h_C(m)\}, \qquad h_F(m)=\frac{mL}{n_F}
\]
where $m$ is a positive integer, and 
{
\begin{align*}
h_C(m) &=\max\{|x_2^{im}-x^{(i+1)m}_2|: i=0,1,\cdots,\lfloor n_C/m\rfloor-1\} \\
    & =\frac{m\pi}{2n_C}\sin(\xi)\;, \ \text{for some} \  \xi \in \left(\frac{2mi-1}{2n_C} \pi, \frac{2m(i+1)-1}{2n_C}\pi\right) \\
    & \approx \frac{m\pi}{2n_C} \\
\end{align*}
where $(x_2^j)$ are the Chebyshev grid points in the $x_2$-direction of the physical space. 
For $n_F\times n_C=256\times128$, this means $h(16)\approx 0.196$ and $h(32) \approx .0.393$. } The nudging parameter is fixed at $\mu=1$.  

\begin{figure}[h] 
	\includegraphics[width=.7\linewidth]{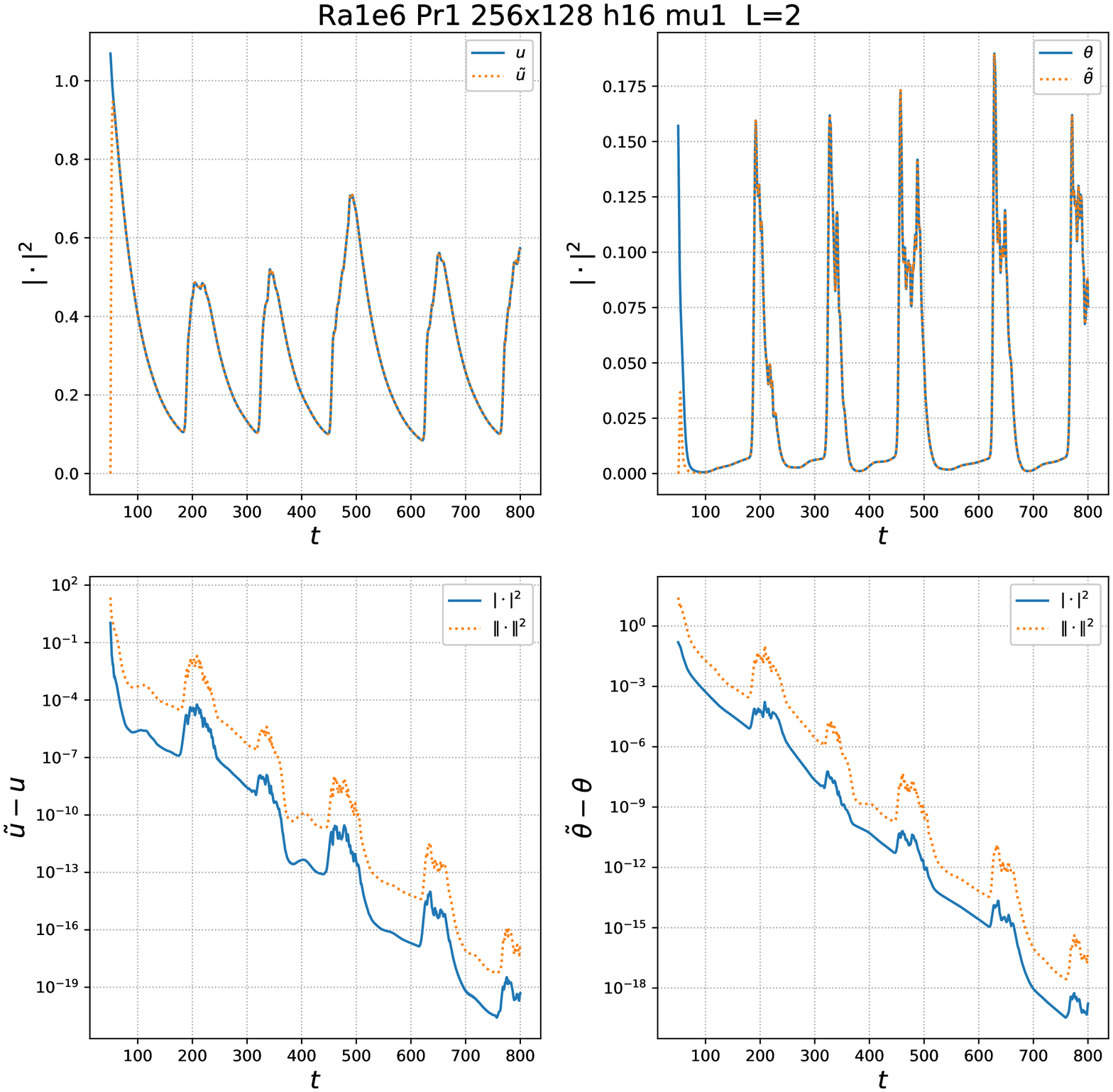}
	\caption{Data assimilation at $\Ra=10^6$ and $\Pr=1$ with $h\approx 0.196$.}
	\label{nudging}
	\bigskip \bigskip 
	\includegraphics[width=.7\linewidth]{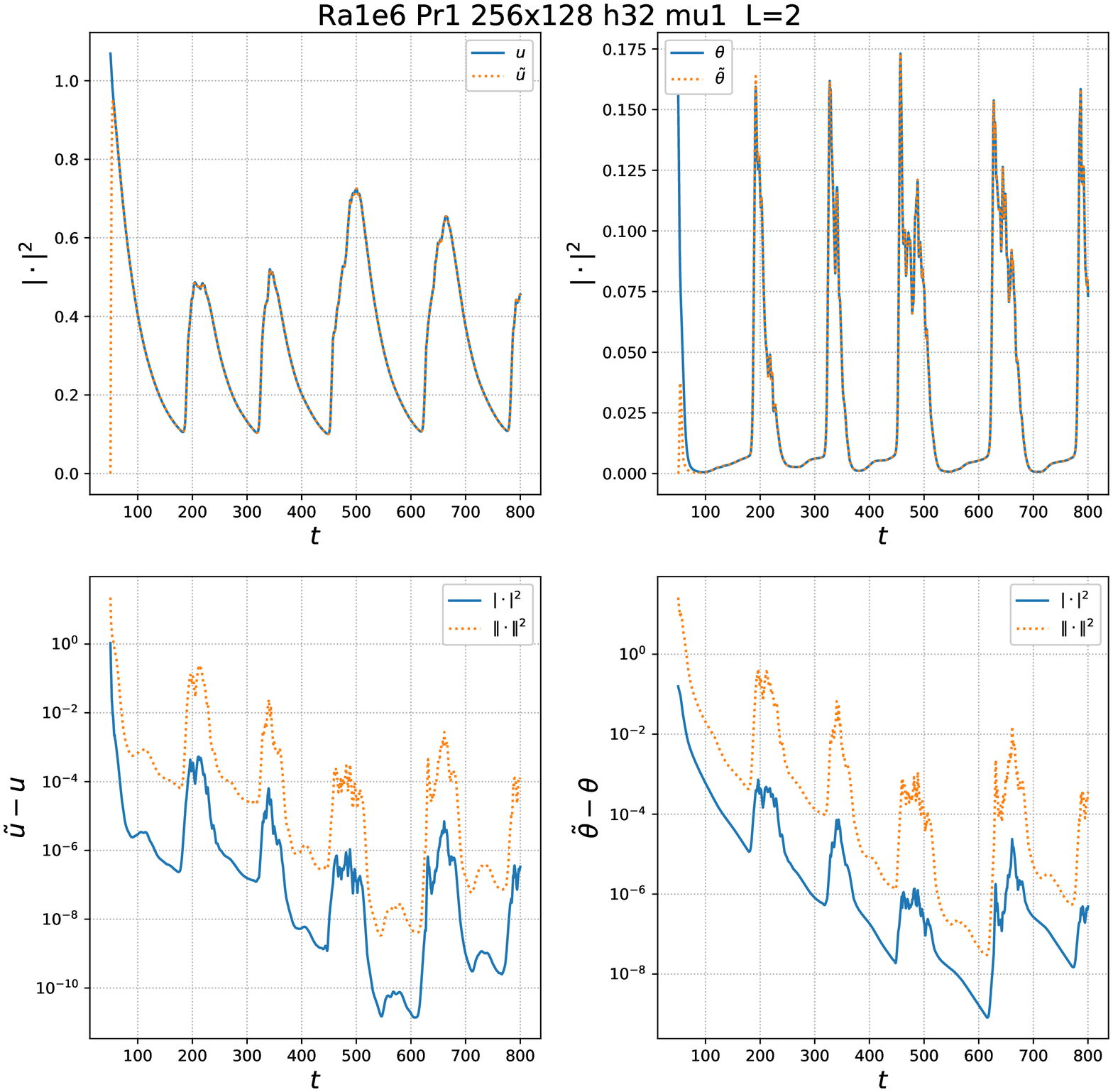}
	\caption{Data assimilation at $\Ra=10^6$ and $\Pr=1$ with $h\approx 0.393$.}
	\label{nudging2}
\end{figure}

Figure \ref{nudging} shows that at $h=0.196$ the solution to the data assimilation system converges to the reference solution at an exponential rate.  At $h=0.393$ the error appears to saturate around $10^{-3}$ during rapid oscillations (see Figure \ref{nudging2}).  We found that at $h=h(64)=0.785$ the nudged solution does not converge to the reference at all (not shown). This demonstrates a critical value of $h$.  

Data assimilation by nudging works much more effectively than the rigorous analysis can guarantee.  The value of $\mu$ and corresponding resolution $h$ of the data suggested by the conditions 
\eqref{condAu} and \eqref{condAth} are based on compounded, conservative estimates derived using general inequalities which are not saturated by 2D convective flows.  In addition, as demonstrated in
\eqref{sharpness}, our algebraic rigorous estimates for $\|\theta\|$, $|A_0u|$, and $|A_1\theta|$ in this case of stress-free boundary conditions, though much better than the exponential bounds previously { found for the no-slip case in \cite{Foias1987attractors}, are still somewhat artificially inflated.}  Numerical nudging tests in \cite{Farhat2018Assim} for the Rayleigh-B\'enard system with no-slip boundary conditions suggest that better bounds on the attractor might hold in that case as well.  { The key here in the stress-free case was extending the physical domain to be fully periodic, hence there is effectively, no boundary.   Since, in the no-slip case one is unable to remove the physical boundary, one should have to resolve the boundary layer scales in order to determine the behavior of the solutions. This is even more pronounced in the estimates of the dimension of the global attractor of the 2D Navier-Stokes equations with no-slip boundary conditions in comparison to the case with periodic boundary conditions. Thus, improving the bounds in the no-slip case would require entirely different techniques.}
\vfill\eject

\vfill\eject

\section{Acknowledgments}   
The authors acknowledge the Indiana University Pervasive Technology Institute (see \cite{stewart2017indiana}) for providing HPC (Big Red II, Carbonate), storage resources that have contributed to the research results reported within this paper.

The work of Y. Cao was supported in part by National Science Foundation grant DMS-1418911, that of M.S. Jolly by  NSF grant DMS-1818754.
The work of E.S. Titi was supported
in part by the Einstein Visiting
Fellow Program, and by the John Simon Guggenheim Memorial Foundation.
J.P. Whitehead acknowledges support from the Simons Foundation through award number 586788, and the hospitality of the Department of Mathematics at Indiana University where part of this work was instigated. 

\bibliographystyle{siam}

\bibliography{CJTWref.bib}

\begin{thebibliography}{10}

\bibitem{Altaf2017Down}
{\sc M.~U. Altaf, E.~S. Titi, T.~Gebrael, O.~M. Knio, L.~Zhao, M.~F. McCabe,
  and I.~Hoteit}, {\em Downscaling the 2{D} {B}\'{e}nard convection equations
  using continuous data assimilation}, Comput. Geosci., 21 (2017),
  pp.~393--410.

\bibitem{Azouani2014continuous}
{\sc A.~Azouani, E.~Olson, and E.~S. Titi}, {\em Continuous data assimilation
  using general interpolant observables}, J. Nonlinear Sci., 24 (2014),
  pp.~277--304.

\bibitem{Azouani2013feedback}
{\sc A.~Azouani and E.~S. Titi}, {\em Feedback control of nonlinear dissipative
  systems by finite determining parameters---a reaction-diffusion paradigm},
  Evol. Equ. Control Theory, 3 (2014), pp.~579--594.

\bibitem{dedalus}
{\sc K.~J. {Burns}, G.~M. {Vasil}, J.~S. {Oishi}, D.~{Lecoanet}, and
  B.~{Brown}}, {\em {Dedalus: Flexible framework for spectrally solving
  differential equations}}.
\newblock Astrophysics Source Code Library, Mar. 2016.

\bibitem{Dascaliuc2010Estimates}
{\sc R.~Dascaliuc, C.~Foias, and M.~S. Jolly}, {\em Estimates on enstrophy,
  palinstrophy, and invariant measures for 2-{D} turbulence}, J. Differential
  Equations, 248 (2010), pp.~792--819.

\bibitem{diLeoni2018Inferring}
{\sc P.~C. Di~Leoni, A.~Mazzino, and L.~Biferale}, {\em Inferring flow
  parameters and turbulent configuration with physics-informed data
  assimilation and spectral nudging}, Phys. Rev. Fluids, 3 (2018), p.~104604.

\bibitem{Farhat2018Assim}
{\sc A.~Farhat, H.~Johnston, M.~Jolly, and E.~S. Titi}, {\em Assimilation of
  nearly turbulent {R}ayleigh-{B}\'{e}nard flow through vorticity or local
  circulation measurements: a computational study}, J. Sci. Comput., 77 (2018),
  pp.~1519--1533.

\bibitem{Farhat2017continuous}
{\sc A.~Farhat, E.~Lunasin, and E.~S. Titi}, {\em Continuous data assimilation
  for a 2{D} {B}\'enard convection system through horizontal velocity
  measurements alone}, J. Nonlinear Sci., 27 (2017), pp.~1065--1087.

\bibitem{Foias2002statistical}
{\sc C.~Foias, M.~S. Jolly, O.~P. Manley, and R.~Rosa}, {\em Statistical
  estimates for the {N}avier-{S}tokes equations and the {K}raichnan theory of
  2-{D} fully developed turbulence}, J. Statist. Phys., 108 (2002),
  pp.~591--645.

\bibitem{Foias1987attractors}
{\sc C.~Foias, O.~Manley, and R.~Temam}, {\em Attractors for the {B}\'enard
  problem: existence and physical bounds on their fractal dimension}, Nonlinear
  Anal., 11 (1987), pp.~939--967.

\bibitem{Gesho2016}
{\sc M.~Gesho, E.~Olson, and E.~S. Titi}, {\em A computational study of a data
  assimilation algorithm for the two-dimensional {N}avier-{S}tokes equations},
  Commun. Comput. Phys., 19 (2016), pp.~1094--1110.

\bibitem{GoJoFlSp2014}
{\sc D.~Goluskin, H.~Johnston, G.~R. Flierl, and E.~A. Spiegel}, {\em
  Convectively driven shear and decreased heat flux}, Journal of Fluid
  Mechanics, 759 (2014), pp.~360--385.

\bibitem{Ibdah2019Uniform}
{\sc H.~A. Ibdah, C.~F. Mondaini, and E.~S. Titi}, {\em {Fully discrete
  numerical schemes of a data assimilation algorithm: uniform-in-time error
  estimates}}, IMA Journal of Numerical Analysis,  (2019).
\newblock https://doi.org/10.1093/imanum/drz043.

\bibitem{Jones1992Determining}
{\sc D.~A. Jones and E.~S. Titi}, {\em Determining finite volume elements for
  the {$2$}{D} {N}avier-{S}tokes equations}, Phys. D, 60 (1992), pp.~165--174.

\bibitem{Mondaini2018Uniform}
{\sc C.~F. Mondaini and E.~S. Titi}, {\em Uniform-in-time error estimates for
  the postprocessing {G}alerkin method applied to a data assimilation
  algorithm}, SIAM J. Numer. Anal., 56 (2018), pp.~78--110.

\bibitem{stewart2017indiana}
{\sc C.~A. Stewart, V.~Welch, B.~Plale, G.~Fox, M.~Pierce, and T.~Sterling},
  {\em Indiana {U}niversity {P}ervasive {T}echnology {I}nstitute},
  https://doi.org/10.5967/K8G44NGB,  (2017).

\bibitem{temam2012infinite}
{\sc R.~Temam}, {\em Infinite-dimensional dynamical systems in mechanics and
  physics}, vol.~68 of Applied Mathematical Sciences, Springer-Verlag, New
  York, second~ed., 1997.

\bibitem{Whitehead2011Ultimate}
{\sc J.~P. Whitehead and C.~R. Doering}, {\em Ultimate state of two-dimensional
  {R}ayleigh-{B}\'enard convection between free-slip fixed-temperature
  boundaries}, Phys. Rev. Lett., 106 (2011), p.~244501.

\end{thebibliography}

\end{document}